\documentclass{article}
\usepackage[utf8]{inputenc}
  \usepackage[active]{srcltx}
\usepackage{amsfonts, amssymb}
\usepackage{amsthm}
\usepackage{amsmath}
\usepackage{enumerate}
\usepackage{graphicx}
\usepackage[labelsep=period,font=small]{caption}
\usepackage{mathrsfs}
\usepackage{authblk}

\def\BB{{\mathcal B}}

\def\PP{{\mathcal P}}

\def\ga{\gamma}

\newcommand{\C}{{{\ensuremath{\mathbb  C}}}}
\newcommand{\Cpct}{\ensuremath{{\operatorname{Cap}}}}

\newcommand{\Chat}{{\ensuremath{\widehat{\mathbb  C}}}}

\newcommand{\Co}{\ensuremath{{\operatorname{Co}}}}

\newcommand{\diam}{\ensuremath{{\operatorname{diam}}}}

\newcommand{\dA}{\ensuremath{d\!A}}

\newcommand{\D}{\ensuremath{{\mathbb   D}}}

\newcommand{\e}{\ensuremath{{\operatorname{e}}}}
\newcommand{\eps}{{\epsilon}}

\newcommand{\Laplace}{{\mathop{{}\bigtriangleup}\nolimits}}

\newcommand{\mapfromto}[3]{\hbox{\ensuremath{#1 \colon #2 \longrightarrow #3}}}

\newcommand{\N}{\ensuremath{\mathbb  N}}

\newcommand{\ooo}{\ensuremath{{\emph{o}}}}

\newcommand{\Reg}{\ensuremath{{\mathcal Reg}}}

\newcommand{\sm}{\ensuremath{{\smallsetminus}}}
\newcommand{\Sm}{\ensuremath{{\setminus}}}

\newcommand{\Supp}{\ensuremath{{\operatorname{Supp}}}}

\def\ga{\gamma}

\def\si{\sigma}

\def\eps{\epsilon}

\def\om{\omega}
\def\Om{\Omega}

\def\C{\mbox{$\mathbb C$}}

\def\N{\mbox{$\mathbb N$}}

\def\D{\mbox{$\mathbb D$}}

\theoremstyle{plain}
\newtheorem{newthm}{Theorem}
\newtheorem{maintheorem}{Theorem}
\newtheorem{notation}{Notation}
\newtheorem{theorem}{Theorem}[section]
\newtheorem{sttheorem}{\cite{StahlandTotik}-Theorem} 
\newtheorem{lemma}[theorem]{Lemma}
\newtheorem{proposition}[theorem]{Proposition}
\newtheorem{corollary}[theorem]{Corollary}
\newtheorem{conjecture}{Conjecture}
\newtheorem{defthm}[theorem]{Definition and Theorem}
\newtheorem{definition}[theorem]{Definition}

\newtheorem{question}[theorem]{Question}

\theoremstyle{remark}
\newtheorem{example}[theorem]{Example}
\newtheorem{remark}[theorem]{Remark}
\newtheorem*{claim1*}{Claim 1}
\newtheorem*{claim2*}{Claim 2}
\newtheorem*{claim3*}{Claim 3}
\newtheorem*{claim4*}{Claim 4}

\newcommand{\ALIGN}{\begin{align*}}
\newcommand{\ENDALIGN}{\end{align*}}
\newcommand{\ENUM}{\begin{enumerate}}
\newcommand{\ENUMa}{\begin{enumerate}[a.]}
\newcommand{\ENUMA}{\begin{enumerate}[A.]}
\newcommand{\ENUMi}{\begin{enumerate}[i)]}
\newcommand{\ENDENUM}{\end{enumerate}}
\newcommand{\ITMZ}{\begin{itemize}}
\newcommand{\ENDITMZ}{\end{itemize}}
\newcommand{\EQN}[1] { \begin{equation}\label{#1} }
\newcommand{\ENDEQN}{\end{equation}}
\newcommand{\THM}{\begin{theorem}}
\newcommand{\EXA}{ \begin{example}}
\newcommand{\REFEXA}[1] { \begin{example}\label{#1} }
\newcommand{\ENDEXA}{\end{example}}
\newcommand{\REM}{ \begin{remark}}
\newcommand{\ENDREM}{\end{remark}}
\newcommand{\REFTHM}[1] { \begin{theorem}\label{#1} }
\newcommand{\RTHM}[1] { \begin{theorem}[#1] }
\newcommand{\RREFTHM}[2] { \begin{theorem}[#1]\label{#2} }
\newcommand{\ENDTHM}{\end{theorem}}
\newcommand{\RREFSTTHM}[2] { \begin{sttheorem}[#1]\label{#2} }
\newcommand{\ENDSTTHM}{\end{sttheorem}}
\newcommand{\REFNTH}[1] { \begin{newthm}\label{#1} }
\newcommand{\ENDNTH}{\end{newthm}}
\newcommand{\REFPROP}[1]{\begin{proposition}\label{#1} }
\newcommand{\RREFPROP}[2]{\begin{proposition}[#1]\label{#2} }
\newcommand{\RPROP}[1]{\begin{proposition}[#1] }
\newcommand{\PROP}{\begin{proposition}}
\newcommand{\ENDPROP}{\end{proposition} }
\newcommand{\REFDEF}[1]{\begin{definition}\label{#1} }
\newcommand{\RREFDEF}[2]{\begin{definition}[#1]\label{#2} }
\newcommand{\DEF}{\begin{definition}}
\newcommand{\RDEF}[1] {\begin{definition}[#1]}
\newcommand{\ENDDEF}{\end{definition} }
\newcommand{\REFLEM}[1]{\begin{lemma}\label{#1} }
\newcommand{\RLEM}[1]{\begin{lemma}[#1] }
\newcommand{\RREFLEM}[2]{\begin{lemma}[#1]\label{#2} }
\newcommand{\LEM}{\begin{lemma}}
\newcommand{\ENDLEM}{\end{lemma} }
\newcommand{\CENTER}{\begin{center}}
\newcommand{\ENDCENTER}{\end{center}}
\newcommand{\REFCOR}[1]{\begin{corollary}\label{#1} }
\newcommand{\RCOR}[1] {\begin{corollary}[#1]}
\newcommand{\COR}{\begin{corollary}}
\newcommand{\ENDCOR}{\end{corollary}}
\newcommand{\REFCONJ}[1]{\begin{conjecture}\label{#1} }
\newcommand{\CONJ}{\begin{conjecture}}
\newcommand{\ENDCONJ}{\end{conjecture}}
\newcommand{\RMRK}{\begin{remark}}
\newcommand{\ENDRMRK}{\end{remark}}
\newcommand{\REFNOTN}[1]{\begin{notation}\label{#1} }
\newcommand{\NOTN}{\begin{notation}}
\newcommand{\ENDNOTN}{\end{notation}}
\newcommand{\REFDEFTHM}[1] { \begin{defthm}\label{#1} }
\newcommand{\RREFDEFTHM}[2] { \begin{defthm}[#1]\label{#2} }
\newcommand{\ENDDEFTHM}{\end{defthm}}
\newcommand{\QUE}{ \begin{question}}
\newcommand{\REFQUE}[1] {\begin{question}\label{#1}}
\newcommand{\ENDQUE}{\end{question}}
\newcommand{\MAINTHM}{\begin{maintheorem}}
\newcommand{\REFMAINTHM}[1] { \begin{maintheorem}\label{#1} }
\newcommand{\ENDMTHM}{\end{maintheorem}}
\newcommand{\ENDCD}{\end{CD}}

\newcommand{\corref}[1]{Corollary~\ref{#1}}

\newcommand{\lemref}[1]{Lemma~\ref{#1}}

\newcommand{\thmref}[1]{Theorem~\ref{#1}}

\newcommand{\propref}[1]{Proposition~\ref{#1}}

\newcommand{\secref}[1]{{Section~\ref{#1}.}}

\newcommand{\PROOF}{\begin{proof}}
\newcommand{\ENDPROOF}{\end{proof}}

\date{\today}

\author{Carsten Lunde Petersen and Eva Uhre}
\affil{Department of Science and Environment, 
Roskilde University, 
4000 Roskilde, 
Denmark
\thanks{
The first author would like to thank 
the Danish Council for Independent Research $|$ Natural Sciences for support via the 
grant DFF -- 4181-00502.}
}
\title{Weak limits of the measures of maximal entropy 
	for Orthogonal polynomials.}
\begin{document}
\maketitle
\begin{abstract}
In this paper we study the sequence of orthonormal polynomials $\{P_n(\mu; z)\}$ 
defined by a Borel probability measure $\mu$ with non-polar compact support $S(\mu)\subset\C$. For each $n\geq 2$ let $\om_n$ denote the unique measure of maximal entropy 
for $P_n(\mu; z)$. We prove that the sequence $\{\om_n\}_n$ is pre-compact for the 
weak-* topology and that for any weak-* limit $\nu$ 
of a convergent sub-sequence $\{\om_{n_k}\}$, the support $S(\nu)$ is contained in 
the filled-in or polynomial-convex hull of the support $S(\mu)$ for $\mu$. 
And for $n$-th root regular measures $\mu$ the full sequence $\{\om_n\}_n$ 
converge weak-* to the equilibrium measure $\om$ on $S(\mu)$.
\end{abstract}

\section{Introduction and general results} 
In the classical study \cite{StahlandTotik} by Stahl and Totik of general orthogonal polynomials 
they relate the potential and measure theoretic properties of the asymptotic zero distribution 
for the sequence of orthonormal polynomials defined by a Borel probability measure 
$\mu$ on $\C$ of infinite, but compact support $S(\mu)$, 
to the potential and measure theoretic properties of $\mu$ and its support.
In the paper \cite{CHPP} Christiansen, Henriksen, Pedersen 
and one of the authors of this paper initiated a study of 
the relation between the potential theoretic properties of $\mu$ 
and the asymptotic (as $n$ tend to $\infty$) potential and measure theoretic properties 
of the Julia sets and filled Julia sets of the orthonormal polynomials $P_n$. 
In this paper we extend this with a study of the weak convergence properties 
of the measures of maximal entropy for the orthonormal polynomials. 

For $\mu$ a Borel probability measure on $\C$ with infinite and compact support $S(\mu)$ 
we denote by $\{P_n(z)\} := \{P_n(\mu; z)\}$ the unique sequence 
\EQN{normalpol}
P_n(z) = \ga_n z^n + \textrm{~lower order terms}.
\ENDEQN
of orthonormal polynomials wrt.~$\mu$.
Then $P_n$ is also characterized as the unique normalized polynomial 
of the form \eqref{normalpol} which is orthogonal to all lower degree polynomials or equivalently for which the monic polynomial $p_n(z) = P_n(z)/\ga_n$ is the unique 
monic degree $n$ polynomial of minimal norm in $L^2(\mu)$.

For $S\subset \C$ a compact non-polar subset such as $S(\mu)$ above 
we denote by $\Om = \Om(S)$ the unbounded connected component of $\C\sm S$, 
by $K = K(S) := \C\sm\Om$ the filled-in or just filled $S$ 
(denoted by the polynomial convex hull in e.g.~\cite{StahlandTotik}), 
by $J = J(S) = \partial K = \partial\Om \subset S$ the outer boundary of $S$, 
by $g_\Om$ the Green's function with pole at infinity for $\Om$ and finally 
by $\om_S$ the equilibrium measure for $S$ (with $\Supp(\om_S) = J$). 
For a measure $\mu\in\BB$ we shall write $K(\mu)$ for $K(S(\mu))$, $J(\mu)$ 
for $J(S(\mu))$ and $\Om(\mu)$ for $\Om(S(\mu)$ or simply $K, J$ and $\Om$, 
when the measure $\mu$ is understood from the context.

\DEF
We denote by $\BB$ the set of Borel probability measure on $\C$ with compact non-polar 
support.

We denote by $\Reg$ the set of $n$-th root regular measures in $\BB$, 
that is
$$
\Reg := \left\{ \mu\in\BB\;\big| 
\lim_{n\to\infty} (\ga_n)^{1/n} = \frac{1}{\Cpct(S(\mu))} \right\}.
$$
\ENDDEF

For $\mu\in\BB$ and $P_n$ the associated sequence of orthonormal polynomials, 
we denote by $\Om_n$ the attracted basin of $\infty$ for $P_n$, 
by $K_n = \C\sm \Om_n$ and $J_n =\partial K_n = \partial\Om_n$ 
respectively the filled Julia set and the Julia set of $P_n$, 
see \secref{Background} for definitions of Julia set and filled Julia set of a polynomial, 
by $g_n$ the Green's function with pole at $\infty$ for $\Om_n$ 
and by $\om_n$ the equilibrium measure for $J_n$ 
or equivalently the measure of maximal entropy for $P_n$ 
(see \cite{Brolin}) and see also below for definitions of these terms). 

Our main result is 
{
\renewcommand*{\themaintheorem}{\Alph{maintheorem}}
\REFMAINTHM{Weakstarlmtsofomn}
Let $\mu\in\BB$, then the sequence $\{\om_n\}_{n\geq 2}$ is pre-compact for 
the weak-* topology and for any limit measure $\nu$ 
for a weakly convergent sub-sequence $\{\om_{n_k}\}_k$ 
$$
S(\nu) \subseteq K(\mu).
$$
Moreover if $\mu\in\Reg$ then
$$\om_{n} \overset{\textrm{weak}^*}\longrightarrow \om_{S(\mu)}.
$$
\ENDMTHM
}
The first part of the theorem is an analogue of the last statement 
in the following Theorem from \cite{StahlandTotik}, 
recast in the above notation.
\RREFTHM{{\cite[Thm.~2.1.1, the first part of which was first proven by Fej{\'e}r \cite{Fejer}]{StahlandTotik}}}{zerocontrol}
Suppose $\mu\in\BB$. All zeros of the orthonormal polynomials $P_n(z)$, $n\in\N$, 
are contained in the convex hull $\Co(K)$, and for any compact subset 
$V\subseteq \Om$ the number of zeros of $P_n(z)$, $n\in\N$, 
on $V$ is bounded as $n\to\infty$. 
Consequently any weak* limit point of (the normalized counting 
meauseres for) the zeros of $P_n$ is supported on $K(\mu)$.
\ENDTHM
The second part of the theorem should be compared to 
\cite[Thm.~3.1.4]{StahlandTotik} and \cite[Thm.~3.6.1]{StahlandTotik}. 
Such a comparison reveals that the 
equilibrium measures $\om_n$ on the Julia sets of $P_n$
have much stronger convergence properties than the counting measure on 
the roots of $P_n$, at least in the case of regular measures $\mu$.

Also the main Theorem is a measure theoretic version of 
the following theorem from \cite{CHPP}, 
which says that under a mild extra condtion (see $\BB_+$ below) 
in capacity any limit of a convergent subsequence $K_{n_k}$ 
(convergent in the Hausdorff topology on the space of compact subsets of $\C$) 
is contained in $K(\mu)$. 
{
\RTHM{{\cite[Thm.1.3]{CHPP}}}
For $\mu\in\BB_+ := \{\mu\in\BB \,|\, \displaystyle{\limsup_{n\to\infty}\; |\ga|_n^{1/n}} < \infty \}$ we have
$$
\limsup_{n\to\infty} K_n \subseteq \Co(K).
$$
Moreover, for any $\eps >0$ and $V_\eps := \{z\in \C \,|\, g_\Om(z) \geq \eps \}$, 
$$
\lim_{n\to\infty}\Cpct(V_\eps\cap K_n) = 0.
$$
\ENDTHM
}

Our proof follows the approach of Stahl and Totik for \thmref{zerocontrol}. 

\section{Background.}\label{Background}
\subsection{Potential theory}
We use the notation of Randsford \cite{Randsford}. 
For $\mu$ a Borel probability measure on $\C$ with compact support 
we define its potential as the sub-harmonic function
$$
p_\mu(z) := \int_{\C}\log|z-w| d\mu(w) = \log|z| + \ooo(1),
$$
which is harmonic on the complement of $S(\mu)$. 
And we define its energy as 
$$
I(\mu) :=  \int_{\C\times\C}\log|z-w| d\mu(w)d\mu(w) = \int_{\C}p_\mu(z) d\mu(z), 
$$
it satisfies $-\infty\leq I(\mu) \leq\log\diam(S(\mu))$.
For a Borel set $B\subset\C$ we denote by $\PP(B)$ the set of probability measures $\mu$ 
with compact support $\Supp(\mu)\subset B$, the energy of $B$ as 
$$
I(B) := \sup\{I(\mu)| \mu\in\PP(B)\}
$$
and its (logarithmic) capacity as 
$$
\Cpct(B) := \e^{I(B)}.
$$
The set $B$ is called polar if $\Cpct(B) = 0$ or equivalently $I(B) = -\infty$ 
and the set $B$ is non-polar otherwise.
For $K\subset\C$ a non-polar compact subset there is a unique measure 
denoted by $\om_K$, which realises the supremum above. 
It is called the equilibrium measure for $K$ and its support is contained 
in the outer boundary of $K$.
And according to Frostman's Theorem, \cite[Thm.~3.3.4]{Randsford} 
its potential $p_\mu$ is bounded from below by $I(K) = I(\om_K)$ and equals 
$I(K)$ on $K$ except for an $F_\si$ polar subset $E$ of $K$. 
The Green's function for $\Om := \C\Sm K$ is the non-negative sub-harmonic 
function 
$$
g_\Om(z) = p_{\om_K}(z)-\log(\Cpct(K)) = p_{\om_K}(z)-I(K).
$$ 
The set $K\Sm E$ is precisely the set of Dirichlet regular boundary points of $K$ 
\cite[4.4.9]{Randsford}. 

\subsection{Polynomial dynamics}
For $P(z) = \ga z^d + \ldots$ a polynomial of degree $d>1$, an easy computation shows 
there exists $R = R_P>0$ such that for any $z$ with $|z|>R$ : $|P(z)| \geq 2|z|$. 
Thus the orbit of such $z$ under iteration converge to $\infty$.
We denote by $\Om_P$ the basin of attraction for $\infty$ for $P$, that is,
\EQN{basinofinfty}
\Om_P := \{z\in\C \,|\, P^k(z) \underset{k\to\infty}\longrightarrow \infty \} 
= \bigcup_{k\geq 0} P^{-k}(\C\sm\overline{D(R)}),
\ENDEQN
where $P^k = \overset{k \textrm{ times}}{\overbrace{P\circ P \circ \ldots \circ P}}$. 
It follows immediately that $\Om_P$ is open and completely invariant, 
i.e.~$P^{-1}(\Om_P) = \Om_P = P(\Om_P)$. 
Denote by $K_P := \C\sm\Om_P\subseteq\overline{D(R)}$ the filled Julia set for $P$ and 
by $J_P := \partial\Om_P = \partial K_P$ the Julia set for $P$. 
Then $K_P$ and $J_P$ are compact and also completely invariant. 
Clearly any periodic point, i.e.~solution of an equation $P^k(z) = z$, $k\in\N$ belongs to $K_P$, 
so that $K_P$ is non-empty. 
It follows from \eqref{basinofinfty} that the filled Julia set $K_P$ can also be described as the nested intersection
\EQN{filledJuliaset}
K_P = \bigcap_{k\geq 0} P^{-k}(\overline{\D(R)}).
\ENDEQN
We denote by {\mapfromto {g_P} \C {[0,\infty)}} 
the Green's function for $\Om_P$ with pole at infinity. 
It follows from \eqref{filledJuliaset} that the Green's function $g_P$ satisfies
$$
g_P(z) = \lim_{k\to\infty} \frac{1}{d^k}\log^+(|P^k(z)|/R) = \lim_{k\to\infty} \frac{1}{d^k}\log^+|P^k(z)|. 
$$
Thus $g_P$  vanishes precisely on $K_P$, i.e.~the exceptional set $E$ for $g_P$ is empty 
and hence every point of $J_P$ is a Dirichlet regular boundary point of $\Om_P$. 
Moreover 
$$g_P(P(z)) = d\cdot g_P(z)
\quad\textrm{and}\quad
\Cpct(K_P) = \frac{1}{|\ga|^{\frac{1}{d-1}}}.
$$
According to Brolin, \cite{Brolin}, the equilibrium measure $\om_P$ for the (filled) Julia set $J_P$ equals the unique measure of maximal entropy for $P$.

\section{Proof of \thmref{Weakstarlmtsofomn}.}
Recall that for $\mu\in\BB$ and $P_n$ the associated orthonormal polynomials, 
we denote by $\Om_n$, $J_n$, and $K_n$, respectively, 
the basin of attraction for $\infty$, the Julia set, and the filled Julia set of $P_n$. 
And moreover by $g_n$ we denote the Green's function for $\Om_n$ 
and by $\om_n$ the equilibrium measure for $K_n$ with support $J_n$. 

The following Lemma gives an elementary relation between the filled support 
$K(\mu)$ and the filled Julia sets $K_n$
\RREFLEM{{\cite[Lemma 2.2]{CHPP}}}{basic_lemma}
Let $\mu\in\BB$ and choose $R>0$ such that $K(\mu) \subset \D(0,R)$. 
Then there exists $N$ such that for all $n\geq N$:
$$
K_n \subset P_n^{-1}(\overline{\D(0, R)})\subset \D(0, R).
$$
\ENDLEM

\REFLEM{equivalentREGconditions}
For a probability measure $\mu\in\BB$ 
the following are equivalent
\ENUM
\item
$$
\lim_{n\to\infty} (\ga_{n})^{1/n} = \frac{1}{\Cpct(S(\mu))}.
$$
\item
$$
\Cpct(K_{n}) \underset{n\to\infty}\longrightarrow \Cpct(S(\mu)).
$$
\item
$$
I(\om_{n}) \underset{n\to\infty}\longrightarrow I(\om_{S(\mu)})
$$
\ENDENUM
\ENDLEM
\PROOF
Since $\Cpct(K_n) = \frac{1}{|\ga|^{1/(n-1)}} = \e^{I(\om_n)}$ 
the lemma follows from the observation that 
$$
\lim_{n\to\infty} (\ga_{n})^{1/{n}} = \frac{1}{\Cpct(S(\mu))} \Longleftrightarrow
\lim_{n\to\infty} (\ga_{n})^{1/(n-1)} = \frac{1}{\Cpct(S(\mu))}.
$$
\ENDPROOF

We shall use the following Lemma which is a refinement of \cite[Lemma 1.3.2]{StahlandTotik} 
\REFLEM{deepbound}
Let $V, S\subset\C$ be two compact subsets with $V\cap K(S) = \emptyset$ and
let $b, 0 < b < 1$ be arbitrary. 
Then there exists $M = M(b, V, S) \in\N$ such that for 
$M$ arbitrary points $x_1, x_2, \ldots x_M \in V$ there exists $M$ points 
$y_1, y_2, \ldots y_M\in \C$ for which the rational function 
$$
r(z) = \prod_{j=1}^M \frac{z-y_j}{z-x_j}
$$
has supremum norm on $S$ bounded by $b$:
$$
||r||_S \leq b.
$$
\ENDLEM
\PROOF 
The above cited Lemma~1.3.2 in \cite{StahlandTotik} shows that the statement holds not for arbitrary $b$, 
but for some value of $b$ call it $a,~0<a<1$ 
and some corresponding value $m = m(V, S)\in\N$. 

To pass from this to the general statement of the lemma 
let $a$ and $m$ be given as in Stahl and Totiks Lemma, 
and let $b, 0<b<1$ be arbitrary. Choose $n\in\N$ such that $a^n<b$ and set $M=nm$. 
Given $x_1, \ldots x_M \in V$ apply \cite[Lemma 1.3.2]{StahlandTotik} to each of the 
$n$ sets of points $\{x_{j+mk}| 1\leq j \leq m\}$, $0\leq k<n$ and multiply. 
\ENDPROOF

\REFPROP{boundonnoofpreimages}
For any $R>0$, for any $\mu\in\BB$ and any compact set $V\subset \C$ with $V\cap K(S(\mu))=\emptyset$ there exists $M=M(R,\mu, V) \in\N$ 
such that for any $w, |w| \leq R$ and 
for any $n$ the number of pre-images of $w$ in $V$ 
under the orthonormal polynomial $P_n$ is less than $M$. 
That is
$$
\#(P_n^{-1}(w) \cap V) < M.
$$
\ENDPROP
\PROOF
Fix $R>0$ and any compact set $V\subset \C$ with $V\cap K(S(\mu))=\emptyset$. 
Let $b = 1/(1+R)$ and let $M=M(b, V, S(\mu))$ be as in \lemref{deepbound} above. 
For each $n$ let $p_n(z) = P_n(z)/\ga_n$ denote the monic orthogonal polynomial 
of degree $n$ for $\mu$. 
Then $p_n$ is the unique degree $n$ monic polynomial 
of minimal norm $1/\ga_n$ in $L^2(\mu)$. 
Suppose towards a contradiction that for some $w, |w|\leq R$ and some $n$ the 
equation $P_n(z)-w$ has at least $M$ solutions $x_1, \ldots x_M \in V$. 
Let $r$ be the rational function given by \lemref{deepbound} such that 
$||r||_{S(\mu)} \leq b$ and set $q(z) := r(z)\cdot(P_n(z)-w)/\ga_n$. 
Then $q$ is a monic polynomial of degree $n$.
Since $P_n$ is orthogonal to all polynomials of lower degree, and in particular to the 
constant polynomial $w$ we compute
\ALIGN
||q||_{L^2(\mu)} &\leq \frac{||r||_{S(\mu)}}{\ga_n} ||P_n(z)-w||_{L^2(\mu)}
\leq\frac{b}{\ga_n}\sqrt{1+|w|^2}\\
&\leq\frac{\sqrt{1+R^2}}{1+R}\cdot ||p_n||_{L^2(\mu)} < ||p_n||_{L^2(\mu)}.
\end{align*}
This contradicts the fact that 
$p_n$ is the degree $n$ monic polynomial with minimal $L^2(\mu)$ norm.
\ENDPROOF

\REFCOR{Knboundedness} 
For any $\mu\in\BB$ and any compact set $V\subset \C$ with $V\cap K(\mu)=\emptyset$ 
there exists $M\in\N$ such that any of the orthonormal polynomials $P_n$, $n\geq 2$ and 
any $z\in K_n$ the number of pre-images of $z$ in $V$ under $P_n$ is less than $M$. 
That is
$$
\#(P_n^{-1}(z) \cap V) < M.
$$
\ENDCOR
\PROOF
It follows from \lemref{basic_lemma} that there exists $R>0$ such that 
$$
K_n \subset P_n^{-1}(\overline{\D(R)}) \subset \D(R),
$$ 
for all $n\geq 2$, where $K_n$ is the filled Julia set of $P_n$. 
As $K_n$ is invariant, $P_n^{-1}(K_n) = K_n$. 
The Corollary follows immediately from this and \propref{boundonnoofpreimages} above.
\ENDPROOF
\PROOF (of \thmref{Weakstarlmtsofomn})
According to \lemref{basic_lemma} the sequence of equillbrium measures $\om_n$ for $K_n$ 
have uniformly bounded support and so the sequence of such measures is pre-compact for 
the weak-* topology.
Furthermore according to Brolin, \cite{Brolin} the equilibriums measure $\om_n$ is also the unique invariant balanced measure for $P_n$,  
i.e.~it is the unique probability measure $\om$ on $\C$ such that for any measurable function 
{\mapfromto f \C \C} 
\EQN{balanced}
\int_{\C} f(z)d\om(z) = \frac{1}{n}\int_{\C} \left(\sum_{w, P_n(w) = z} f(w)\right) d\om(z).
\ENDEQN
Let $V\subset\C$ be a compact subset with $V\cap K(\mu) = \emptyset$ 
and let $M\in\N$ be as in \corref{Knboundedness}. 
Then for the measureable function $1_V$, the indicator function for $V$ we have 
$$
\om_n(V) = \int_ {\C}1_V(z)d\om_n(z) = \frac{1}{n}\int_{\C} \left(\sum_{w, P_n(w) = z} 1_V(w)\right) d\om_n(z) 
\leq \frac{M}{n}\underset{n\to\infty}\longrightarrow 0.
$$
This proves the first statement, that for any weak limit $\nu$ 
of a convergent subsequence $\{\om_{n_k}\}_k$ : $S(\nu) \subseteq K(\mu)$.

Next suppose that $\mu\in\Reg$. 
Then in particular $\mu\in\BB_+$, i.e. $ \displaystyle{\limsup_{n\to\infty}\; |\ga|_n^{1/n}} < \infty$ and so according to 
\cite[Cor. 4.1 \textit{iii)}]{CHPP} 
$$
\limsup_{n\to\infty} g_n(z) \leq g_\Om(z).
$$
where $g_n$ is the Green's function with pole at infinity for $\Om_n = \C\Sm K_n$ and $g_\Om$ 
is the Green's function with pole at infinity for $\Om = \C\Sm K(\mu)$. 
Since $g_\Om\equiv 0$ on the interior of $K(\mu)$ we also have 
$\lim_{n\to\infty} g_n \equiv 0$ uniformly (in the usual sense) 
on any compact subset $L$ of the interior of $K(\mu)$. 

Since $\om_n$ can be retrieved from the LaPlacian of $g_n$ in distributions, it follows 
that $\Supp(\nu)\subset J = \partial\Om$. In more detail. 
For any at least $C^2$ test-function $\phi$ with compact support $L$ 
contained in the interior of $K(\mu)$ 
we have
$$
\int_{\C} \phi(z)d\om_n(z) = \frac{1}{2\pi}\int_{\C}\Laplace\phi(z) g_n(z)\dA(z) 
\underset{n\to\infty}\longrightarrow 0
$$
where $\dA(z)$ is the standard Euclidean area element. Thus 
$$
\int_{\C}\phi(z)d\nu(z) = \lim_{n\to\infty}\int_{\C} \phi(z)d\om_n(z) = 0.
$$
It follows that $\nu(L) = 0$ for any compact subset $L$ 
of the interior $\overset\circ{K}(\mu)$ of $K(\mu)$. 
Take e.g.~a compact neighbourhood $M$ of $L$ contained in the interior of $K(\mu)$ 
and let $\phi$ be an at least $C^2$ function with support in $M$, with $0\leq \phi \leq 1$ 
and with $\phi\equiv 1$ on $L$. Then $\nu(L)\leq \int_{\C}\phi(z)d\nu(z) = 0$.
Finally 
$$
\nu(\overset\circ{K}(\mu) = \sup\{\nu(L)| L\subset\overset\circ{K}(\mu), 
L\textrm{~compact}\} = 0
$$ by regularity of $\nu$ as a Borel, in fact Radon measure. 

Finally we have, by the definition of $n$-th-root regularity of $\mu$, \lemref{equivalentREGconditions} and \cite[Lemma 3.3.3]{Randsford})
$$
I(\nu) \geq \limsup_{n\to\infty} I(\om_n) = I(\om_{S(\mu)})
$$
Hence $\nu=\om_{S(\mu)}$, since $\om_{S(\mu)}$ is the unique measure of maximal 
energy $I(S(\mu))$.\\
Since there is a unique possible limit we have in fact 
$$\om_n \overset{\textrm{weak}^*}\longrightarrow \om_{S(\mu)}.$$
\ENDPROOF

\end{document}